\documentclass[aip,
amsmath,amssymb,numerical,preprint,%
]{revtex4-1}
\usepackage{graphicx}
\usepackage{dcolumn}
\usepackage{bm}
\usepackage{hyperref}
\usepackage[utf8]{inputenc}
\usepackage[T1]{fontenc}
\usepackage{mathptmx}
\usepackage{etoolbox}
\newcommand{\RNum}[1]{\uppercase\expandafter{\romannumeral #1\relax}}
\makeatletter
\def\@email#1#2{%
	\endgroup
	\patchcmd{\titleblock@produce}
	{\frontmatter@RRAPformat}
	{\frontmatter@RRAPformat{\produce@RRAP{*#1\href{mailto:#2}{#2}}}\frontmatter@RRAPformat}
	{}{}
}%
\makeatother

\begin{document}
	
	\title{Phototactic bioconvection in an algal suspension with a free top wall due to diffuse flux in the absence of direct collimated flux}
	\author{S. K. Rajput}
	\altaffiliation[Corresponding author: E-mail: ]{shubh.iiitj@gmail.com.}
	\author{M. K. Panda}%
	\affiliation{ 
		Department of Mathematics, PDPM Indian Institute of Information Technology Design and Manufacturing,
		Jabalpur 482005, India.
	}%
	
	
	
	
	\begin{abstract}
		In this article, the effect of diffuse flux in the absence of direct collimated flux on the onset of phototactic bioconvection is investigated. The main effect of diffuse flux in the absence of collimated flux is on the swimming behaviour of microorganisms in the suspension. At higher diffuse flux, the horizontal component of swimming orientation exhibits a higher magnitude, which slows the rate of pattern formation in the suspension. Also, the linear stability of the suspension predicts that the most unstable mode of disturbance transits from stationary (oscillatory) to oscillatory (stationary) at the variation in the magnitude of diffuse flux for some fixed parameters. The overstable nature of disturbance is mostly observed at the high value of swimming speed and extinction coefficient. Moreover, the suspension shows a more stable behaviour at the higher magnitude of diffuse flux.
		
	\end{abstract}
	
	
	\maketitle
	
	
	\section{INTRODUCTION}
	The bioconvection (or biological convection) is a fascinating phenomenon and it occurs due to the macroscopic convective motion of a fluid containing self-propelled motile microorganisms (algae and bacteria)~\cite{20platt1961,21pedley1992,22hill2005,23bees2020,24javadi2020}. In most cases, microorganisms are slightly denser than their surrounding medium (usually water) and swim upwards on average. However, up-swimming and higher density are not necessary for the development of biological patterns~\cite{21pedley1992}. Also, when the microorganisms stop swimming the biological pattern disappears. The swimming microorganisms change their swimming direction in response to different types of stimuli. These responses are known as $taxes$~\cite{7ghorai2010,15panda2016}. $Gravitaxis$ is the swimming response due to gravitational acceleration; $chemotaxis$ is the swimming response due to chemicals; $gyrotaxis$ is caused by the balance between a couple of torques due to gravity and local shear flow; and  $phototaxis$ is defined as the swimming response to the light source.
	
	Direct (vertical and oblique) and indirect (diffuse/scattered) light, can impact the wavelength of the pattern formed in bioconvection~\cite{1wager1911,2kitsunezaki2007,3kessler1985,4williams2011,27kage2013,25kessler1986,26kessler1997,28mendelson1998}. If the light intensity is high, it dissipates steady patterns or intercepts pattern formation in well-stirred cultures. The intensity of light can also alter the shape and size of bioconvection patterns~\cite{3kessler1985,4williams2011,5kessler1989,37ghorai2009}. These changes can be due to several reasons, which can be described as follows. First, motile microorganisms such as algae acquire energy via the process of photosynthesis and move towards the region of strong light intensity (positive phototaxis) when $G< G_c$ and towards the region of low light intensity (negative phototaxis) when $G> G_c$ to avoid photo-damage. This behaviour leads to the accumulation of cells at a suitable place in their natural environment where the light intensity $G\approx G_c$~\cite{6hader1987}. Second, the absorption and scattering can cause changes in pattern formation~\cite{7ghorai2010}. Also, diffuse flux can affect the formation of patterns due to the uniformity of light in the suspension.
	
	In this study, we employ a phototaxis model developed by Panda $et$ $al.$~\cite{15panda2016}. Diffuse light is a kind of sunlight created by the scattering of direct sunlight by water droplets (clouds) and dust particles in the atmosphere. In some atmospheric conditions, only a diffuse part of sunlight can reach the earth's surface. So, the importance of diffuse irradiation should not be ignored in the process of photo-bioconvection. Photo-bioconvection plays an important role in carbon dioxide fixation and biofuel production in photo-bioreactors~\cite{30khan2017,31hayat2021}. Therefore, when designing efficient photo-bioreactors, the effect of diffuse flux in the absence of direct (collimated) flux on phototaxis and the resulting bioconvection must be taken. Due to the need for photosynthesis, the algae are primarily phototactic in their natural environment. Therefore, a realistic phototaxis model needs to account for the effects of diffuse flux in the absence of direct (collimated) flux to fully understand their swimming behaviour.
	
	Consider bioconvection in dilute phototactic algal suspension at a small depth $H$. The basic state (sublayer) of the same suspension is formed, where the fluid becomes motionless and the cell's motion is driven by the interplay of phototaxis and diffusion of cells. The sublayer formation in the suspension depends on the critical intensity $G_c$ and can occur at different depths. If $G<G_c$, the sublayer forms at the top, whereas if $G>G_c$, it forms at the bottom of the suspension. When $G=G_c$ occurs at a fixed depth inside the suspension, a sublayer forms in between the suspension and divides the suspension into two subregions, with the upper region being gravitationally stable and the lower region gravitationally unstable. The sublayer's instability leads to fluid flow in the lower region, which infiltrates the stable upper region, a phenomenon known as penetrative or infiltrative convection. This phenomenon has been observed in a variety of several other convection problems~\cite{9straughan1993,10ghorai2005,11panda2016}.
	
	\begin{figure}[!h]
		\centering
		\includegraphics[width=12.5cm]{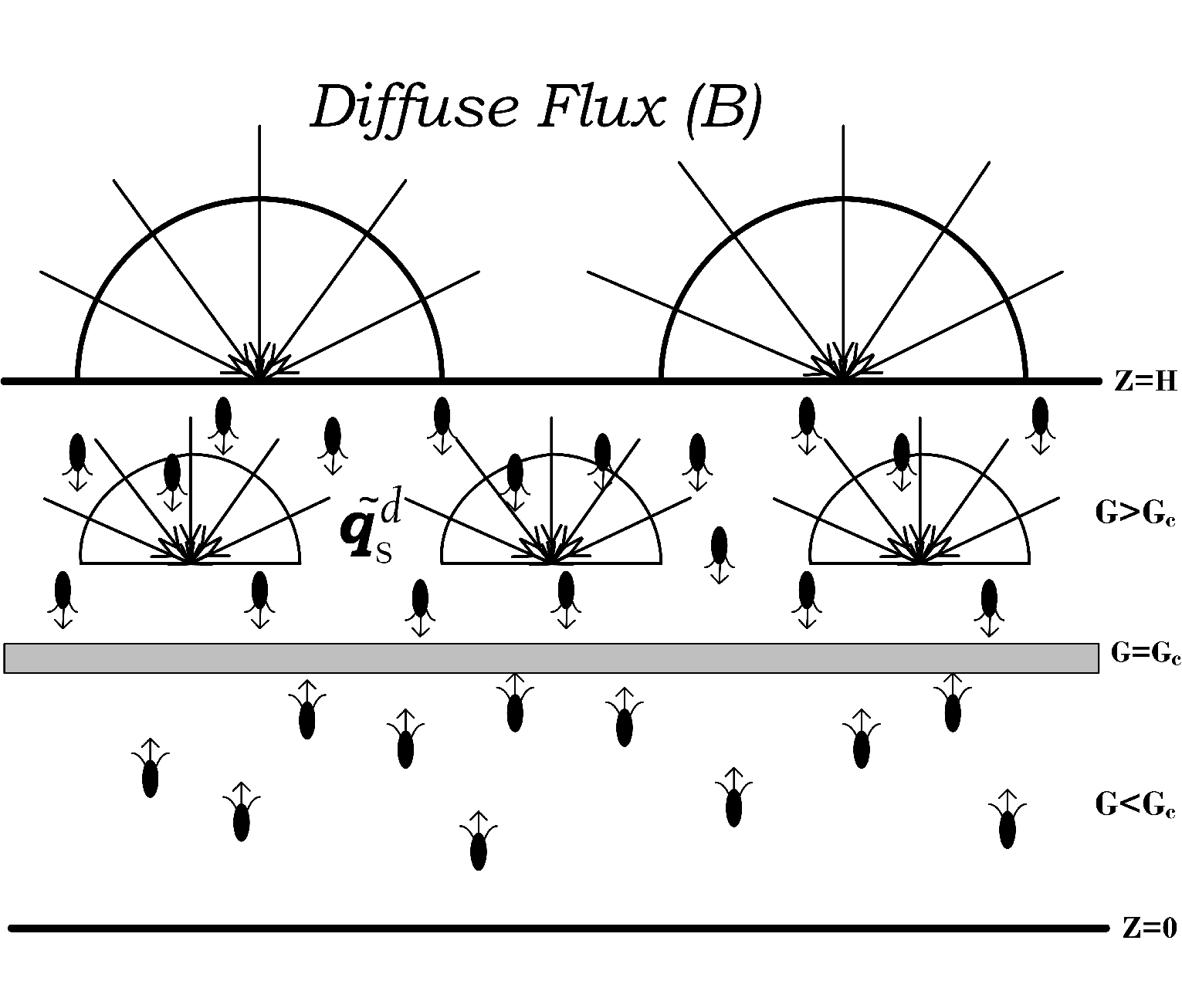}
		\caption{\footnotesize{Formation of sublayer inside the suspension at a depth where $G=G_c$. }}
		\label{fig1}
	\end{figure}
	
	The phototactic bioconvection has been investigated by several researchers in different conditions. The onset investigation in this field was carried out by Vincent and Hill~\cite{12vincent1996}. By using linear stability analysis, they found stationary as well as oscillatory bio-convective solutions at the initiation of convective instability. After that, Ghorai and Hill~\cite{10ghorai2005} investigated the behaviour of phototactic algal suspension in two dimensions but did not consider the scattering effect. Ghorai $et$ $al$.~\cite{7ghorai2010} developed a model in which they considered the scattering effect isotropically and found a bimodal steady-state for highly scattering algal suspension and oscillatory instabilities for certain parameters. Ghorai and Panda~\cite{13ghorai2013} studied the impact of anisotropic scattering in a suspension of phototactic algae and observed that as the variations are made in the anisotropic coefficient, the nature of disturbance change from stationary to oscillatory for some parametric values. Panda and Ghorai~\cite{14panda2013} proposed a model for isotropically scattering medium in two-dimension and found results that were different from those found by Ghorai and Hill~\cite{10ghorai2005} due to scattering effects. Panda and Singh~\cite{11panda2016} investigated phototactic bioconvection in two dimensions in which the suspension was confined between rigid side walls, observing a considerable stabilizing effect on the suspension. Panda $et$ $al$.~\cite{15panda2016} proposed a model to investigate the effect of diffuse flux in the presence of collimated flux. Their study showed that diffuse flux had a stabilizing impact on the system, but the cell swimming orientation is more influenced by collimated irradiation. The impact of anisotropic scattering on the bioconvective instability in the algal suspension in the presence of both diffuse and non-diffuse (direct) collimated irradiation, was studied by Panda~\cite{8panda2020}. Panda $et$ $al$.~\cite{16panda2022}, in their study, analysed the impact of oblique collimated flux and checked the linear stability of the suspension for different angles of incidence. They found the angle of incidence has a significant stabilizing impact on the suspension stability. Kumar~\cite{17kumar2022} explored the impact of collimated oblique flux in scattering algal suspension isotropically and found both types of solutions (stationary and overstable) for certain ranges of parameters. Recently, Panda and Rajput~\cite{41rajput2023}, investigated the concatenated effect of diffuse and oblique flux at the initiation of instability in the algal suspension.  However, no research has yet examined how diffuse flux in the absence of direct flux might trigger bioconvection in the algal suspension. Therefore, this study investigates the effect of diffuse flux without collimated flux on the algal suspension to gain further insight into phototactic bioconvection.
	
	The article is structured as follows: First, the problem is formulated mathematically and then a basic (equilibrium) solution is derived. After that, a small disturbance is made to the equilibrium system, and then the linear stability problem is solved numerically by using the linear perturbation theory. Then the results of the model are shown. The final discussion centres around the outcomes of this model.  
	
	
	\section{MATHEMATICAL FORMULATION OF THE PROBLEM}
	
	Considered a layer of phototactic microorganisms that extends infinitely in the horizontal direction and is bounded by horizontal boundaries located at $z=0$ and $H$. The system is exposed to diffuse irradiation from the top, as depicted in Fig.~\ref{fig1}. The intensity of light at a particular location $\boldsymbol{x}$ is denoted by $L(\boldsymbol{x},\boldsymbol{s})$ in a specific direction $\boldsymbol{s}$, where $\boldsymbol{s}$ is represented by the unit vector $\boldsymbol{s}=\sin\vartheta\cos\varphi\hat{i}+\sin\vartheta \sin\varphi\hat{j}+\cos\vartheta\hat{k}$~\cite{7ghorai2010}. Here, $\vartheta$ and $\varphi$ are the azimuthal and polar angles.
	
	
	\subsection{\label{sec:level3}The mean swimming orientation}
	The distribution of light intensity within suspension that both scatters and absorbs light is determined through the application of the radiative transfer equation (RTE), represented as  
	\begin{equation}\label{1}
	\boldsymbol{s}\cdot\nabla L(\boldsymbol{x},\boldsymbol{s})+(\varkappa+\sigma)L(\boldsymbol{x},\boldsymbol{s})=\frac{\sigma}{4\pi}\int_{4\pi}L(\boldsymbol{x},\boldsymbol{s'})\mathbf{g}(\boldsymbol{s},\boldsymbol{s'})d\Omega',
	\end{equation}
	where, $\varkappa$ and $\sigma$ are symbols used to denote the absorption and scattering coefficients, respectively. The function $\mathbf{g}(\boldsymbol{s},\boldsymbol{s'})$ provides information about how light intensity is scattered from direction $\boldsymbol{s'}$ into direction $\boldsymbol{s}$, including its angular distribution. Here, the light is assumed to be scattered isotropically by the algal cells. So, for simplicity, we use $\mathbf{g}(\boldsymbol{s},\boldsymbol{s'})=1$.~\cite{15panda2016}

	%
	Assuming a diffusive boundary condition at the top of the suspension, the intensity at a point $x_H=(x,y,H)$ on the upper boundary can be expressed as 
	\begin{equation*}
		L(\boldsymbol{x}_H,\boldsymbol{s})=\frac{B}{\pi}, 
	\end{equation*}
	where $B$ is the magnitude of the diffuse flux~\cite{8panda2020,15panda2016,41rajput2023}. Now, we assume that  $\varkappa={\varkappa_H}n(\boldsymbol{x})$ and $\sigma={\sigma_s}n(\boldsymbol{x})$, then scattering albedo is defined as $\omega=\sigma_s/\kappa$, where $\kappa=\varkappa_H+\sigma_s$. Therefore, RTE becomes in terms of scattering albedo
	\begin{equation}\label{2}
	\boldsymbol{s}\cdot\nabla L(\boldsymbol{x},\boldsymbol{s})+\kappa n(\boldsymbol{x})L(\boldsymbol{x},\boldsymbol{s})=\frac{\omega \kappa n(\boldsymbol{x})}{4\pi}\int_{4\pi}L(\boldsymbol{x},\boldsymbol{s}')d\Omega'.
	\end{equation}
	At a fixed location $\boldsymbol{x}$ within the medium, the total intensity of light is given by
	\begin{equation*}
		G(\boldsymbol{x})=\int_{4\pi}L(\boldsymbol{x},\boldsymbol{s})d\Omega,
	\end{equation*}
	and the radiative heat flux (hereinafter referred to as RHF) is determined by
	\begin{equation}\label{3}
	\tilde{\boldsymbol{q}}(\boldsymbol{x})=\int_{4\pi}L(\boldsymbol{x},\boldsymbol{s})\boldsymbol{s}d\Omega.
	\end{equation}
	The mean cell
	swimming velocity is defined as~\cite{18hill1997}
	\begin{equation*}
		\boldsymbol{U}_c=U_c<\boldsymbol{p}>,
	\end{equation*}
	where, $U_c$ refers to the mean cell's swimming speed, and $<\boldsymbol{p}>$ represents the mean direction of their swimming, which can be determined  by
	\begin{equation}\label{4}
	<\boldsymbol{p}>=-\frac{\boldsymbol{\tilde{q}}}{|\boldsymbol{\tilde{q}}|}M(G),
	\end{equation}
	where $M(G)$ is the taxis function. The function $M(G)$ is non negative (negative) when $G(\boldsymbol{x})\leq G_{c}$ ($G(\boldsymbol{x})> G_{c}$). The specific mathematical expression of $M(G)$ can vary based on the species of microorganism considered~\cite{12vincent1996}. For a particular instance, $M(G)$ can be expressed mathematically as 
	\begin{equation}\label{5}
	M(G)=0.8\sin\bigg[\frac{3\pi}{2}\aleph(G)\bigg]-0.1\sin\bigg[\frac{\pi}{2}\aleph(G)\bigg],
	\end{equation}
	where $\aleph(G)=\frac{G}{2.5}\exp[\varpi(G-2.5)]$ and $\varpi$ depends on critical light intensity.
	
	
	\subsection{The governing equations}
	Similar to previous continuum models~\cite{21pedley1992,10ghorai2005,11panda2016,12vincent1996,14panda2013,15panda2016,16panda2022,17kumar2022,22hill2005,23bees2020,3kessler1985,24javadi2020,1wager1911,2kitsunezaki2007,4williams2011,5kessler1989,7ghorai2010,9straughan1993,8panda2020,41rajput2023}, the suspension is composed of cells with a volume of $v$ and density of $\rho+\delta\rho$ ($\rho$ is water density and $0<\delta\rho\ll\rho$. The velocity and concentration of the cells in the suspension are denoted by $U$ and $n$, respectively. Additionally, it is assumed that the suspension is incompressible, and thus the equation of continuity can be expressed as follows
	\begin{equation}\label{6}
	\boldsymbol{\nabla}\cdot \boldsymbol{U}=0.
	\end{equation}
	The Boussinesq approximation yields the following expression for the momentum equation
	
	\begin{equation}\label{7}
	\rho\left[\frac{DU}{Dt}\right]=-\boldsymbol{\nabla} P+\mu\nabla^2\boldsymbol{U}-nvg\Delta\rho\hat{\boldsymbol{z}},
	\end{equation}
	where $P$ and $\mu$ represent the excess pressure and viscosity of the suspension (or water), respectively.
	
	The equation of cell conservation is given by
	\begin{equation}\label{8}
	\frac{\partial n}{\partial t}=-\boldsymbol{\nabla}\cdot\boldsymbol{F}_0,
	\end{equation}
	where ${\boldsymbol{F}_0}$ is the total cell flux  
	\begin{equation}\label{9}
	\boldsymbol{F}_0=n[U+U_c<\boldsymbol{p}>]-\boldsymbol{D}\boldsymbol{\nabla}n.
	\end{equation}
	The total cell flux consists of two parts, one due to advection and the other due to diffusion. Here, $\boldsymbol{D}$ is the diffusion tensor, which is considered to be constant and isotropic. Hence, $\boldsymbol{D} = DI$. It is also assumed that the cells exhibit purely phototactic behaviour, so the effect of viscous torque can be neglected. These simplifications allow for the removal of the Fokker-Planck equation from the system of governing equations.~\cite{12vincent1996}
	
	
	\subsection{The boundary conditions}
	In this article, it is assumed that the top boundary is stress-free, while the bottom boundary is rigid non-slip. Hence, it is emphasized that there is no flow of fluid and algae cells across the boundaries. Hence, boundary conditions can be formulated as
	\begin{equation}\label{10}
	\boldsymbol{U}\cdot\hat{\boldsymbol{z}}=0,~~~~~z=0,H,
	\end{equation}
	\begin{equation}\label{11}
	\boldsymbol{F}_0\cdot\hat{\boldsymbol{z}}=0,~~~~~z=0,H.
	\end{equation}
	For rigid boundary
	\begin{equation}\label{12}
	\boldsymbol{U}\times\hat{\boldsymbol{z}}=0,~~~~~z=0,
	\end{equation}
	and for stress-free boundary
	\begin{equation}\label{13}
	\frac{\partial^2}{\partial z^2}(\boldsymbol{U}\cdot\hat{\boldsymbol{z}})=0,~~~~~z=H.
	\end{equation}
	Now, let's consider a scenario where the top boundary surface is subjected to a source of diffuse flux uniformly, in the absence of direct collimated flux. Under these conditions, the boundary conditions for intensities can be described as
	\begin{subequations}
		\begin{equation}\label{14a}
		L(x, y, z=H, \vartheta, \varphi) =\frac{B}{\pi},~~where~~ (\pi/2\leq\vartheta\leq\pi),
		\end{equation}
		\begin{equation}\label{14b}
		L(x, y, z=0, \vartheta, \varphi) =0,~~where~~ (0\leq\vartheta\leq\pi/2).
		\end{equation}
	\end{subequations}

	Now, there is a need to scaled the governing equations to make the governing equations dimensionless. Here, we use a set of characteristic scales for length, time, velocity, pressure, and concentration similar to the previous models~\cite{21pedley1992,10ghorai2005,11panda2016,12vincent1996,14panda2013,15panda2016,16panda2022,17kumar2022,22hill2005,23bees2020,3kessler1985,24javadi2020,1wager1911,2kitsunezaki2007,4williams2011,5kessler1989,7ghorai2010,9straughan1993,8panda2020,41rajput2023}. Therefore, the scaled governing equations are expressed as
	\begin{equation}\label{15}
	\boldsymbol{\nabla}\cdot\boldsymbol{U}=0,
	\end{equation}
	\begin{equation}\label{16}
	S_{c}^{-1}\left[\frac{D\boldsymbol{U}}{Dt}\right]=-\nabla P-Rn\hat{\boldsymbol{z}}+\nabla^{2}\boldsymbol{U},
	\end{equation}
	\begin{equation}\label{17}
	\frac{\partial{n}}{\partial{t}}=-{\boldsymbol{\nabla}}\cdot\boldsymbol{F}_0,
	\end{equation}
	where
	\begin{equation}\label{18}
	\boldsymbol{F}_0=n({\boldsymbol{U}}+V_{c}<{\boldsymbol{p}}>)-{\boldsymbol{\nabla}}n,
	\end{equation}

	in which, $S_{c}=\nu/D$ is the Schmidt number, $V_c=U_cH/D$ is the scaled cell swimming speed and $R=\bar{n}vg\Delta{\rho}H^{3}/\nu\rho{D}$ is the Rayleigh number( a control parameter in the bioconvection). In the form of non-dimensional variables, the boundary conditions become
	\begin{equation}\label{19}
	\boldsymbol{U}\cdot\hat{\boldsymbol{z}}=0,~~~~~z=0,1,
	\end{equation}
	\begin{equation}\label{20}
	\boldsymbol{F}_0\cdot\hat{\boldsymbol{z}}=0,~~~~~z=0,1.
	\end{equation}
	For a rigid boundary,	
	\begin{equation}\label{21}
	\boldsymbol{U}\times\hat{\boldsymbol{z}}=0,~~~~~z=0,
	\end{equation}
	and for a stress-free boundary
	\begin{equation}\label{22}
	\frac{\partial^2}{\partial z^2}(\boldsymbol{U}\cdot\hat{\boldsymbol{z}})=0,~~~~~z=1.
	\end{equation}	
	
	In dimensionless form RTE becomes
	\begin{equation}\label{23}
	\boldsymbol{s}\cdot\nabla L(\boldsymbol{x},\boldsymbol{s})+\kappa_Hn(\boldsymbol{x})L(\boldsymbol{x},\boldsymbol{s})=\frac{\omega \kappa_H n(\boldsymbol{x})}{4\pi}\int_{4\pi}L(\boldsymbol{x},\boldsymbol{s}')d\Omega'.
	\end{equation}
	where $\kappa_H=\kappa\Bar{n}H$ denotes the non-dimensional absorption coefficient. In the form of direction cosines $(\alpha,\beta,\gamma)$ of the direction $\boldsymbol{s}$, where
	\begin{equation}
	\alpha=\sin\vartheta\cos\varphi,~~\beta=\sin\vartheta\sin\varphi,~~\gamma=\cos\vartheta
	\end{equation}
	
	RTE becomes,
	\begin{equation}\label{25}
	(\alpha,\beta,\gamma)\cdot\nabla L(\boldsymbol{x},\boldsymbol{s})+\kappa_H nL(\boldsymbol{x},\boldsymbol{s})=\frac{\omega\kappa_Hn }{4\pi}\int_{4\pi}L(\boldsymbol{x},\boldsymbol{s}')d\Omega',
	\end{equation}
	and, the intensity at boundaries is converted to dimensionless form as follows
	
	\begin{equation}\label{26}
	L(x, y, z=1, \vartheta, \varphi) =\frac{B}{\pi},~~where~~ (\pi/2\leq\vartheta\leq\pi),
	\end{equation}

	
	\section{THE STEADY SOLUTION}
	Eqs. $(\ref{15})$-$(\ref{18})$ and $(\ref{25})$, along with appropriate boundary conditions, have a steady state solution in the form of
	\begin{equation}\label{27}
	\boldsymbol{U}=0,\quad n=n_s(z),\quad and\quad  L=L_s(z,\vartheta).
	\end{equation}

	Therefore, the total intensity ($G_s$) and the RHF at the steady state,  $\tilde{\boldsymbol{q}}_s$, can be found by the following equations	
	
	\begin{equation*}
		G_s=\int_{4\pi}L_s(z,\vartheta)d\Omega,\quad 
		{\tilde{\boldsymbol{q}}}_s=\int_{4\pi}L_s(z,\vartheta)\boldsymbol{s} d\Omega,
	\end{equation*}
	and $L_s$ can be governed by the equation
	\begin{equation}\label{28}
	\frac{dL_s}{dz}+\left[\frac{\kappa_H n_s}{\gamma}\right]L_s=\frac{\omega\kappa_H n_s}{4\pi\gamma}G_s(z).
	\end{equation}
	In the steady state, the light intensity consists solely of a diffuse component denoted as $L_s=L_s^d$, and it is calculated by the following equation    
	\begin{equation}\label{29}
	\frac{dL_s^d}{dz}+\left[\frac{\kappa_H n_s}{\gamma}\right]L_s^d=\frac{\omega\kappa_H n_s}{4\pi\gamma}G_s(z),
	\end{equation}
	with the boundary conditions
	\begin{subequations}
		\begin{equation}\label{30a}
		L_s^d( z=1, \vartheta) =\frac{B}{\pi},~~where~~ (\pi/2\leq\vartheta\leq\pi), 
		\end{equation}
		\begin{equation}\label{30b}
		L_s^d( z=0, \vartheta) =0,~~where~~(0\leq\vartheta\leq\pi/2). 
		\end{equation}
	\end{subequations}
	Hence, $G_s$ can be determined using the following relationship
	
	\begin{equation*}
		G_s=G_s^d=\int_{4\pi}L_s^d(z,\vartheta)d\Omega.
	\end{equation*}
	
	Now define a vertical optical depth as
	\begin{equation*}
		\tau(z)=\int_z^1 \kappa_H n_s(z')dz',
	\end{equation*}
	
	Therefore, $G_s$ depends on $\tau$ only. Additionally, non-dimensional total intensity, denoted as $\Upsilon(\tau)=G_s(\tau)$, can be described by the following Fredholm integral equation (FIE)
	\begin{equation}\label{31}
	\Upsilon(\tau)=\frac{\omega}{2}\int_0^{\kappa_H} E_1(|\tau-\tau'|)\Upsilon(\tau')d\tau'+2BE_2(\tau),
	\end{equation}
	where $E_1(x)$ and $E_2(x)$ refer to the exponential integral of the first and second order, respectively~\cite{36chandrasekhar1960}.
	The given FIE exhibits a singularity at $\tau=\tau'$~\cite{35modest2021}, and it is resolved by utilizing the subtraction of singularity method (for detail see Crossbie and Pattabongse~\cite{42crosbie1985,38press1992}).\par
	
	The RHF in the steady state can be expressed as
	\begin{equation*}
		\tilde{{\boldsymbol{q}}}_s=\tilde{{\boldsymbol{q}}}_s^d= \int_{4\pi}L_s^d(z,\vartheta)\boldsymbol{s}d\Omega.
	\end{equation*}
	As $L_s^d$ does not depend on $\varphi$, the horizontal components of $\tilde{\boldsymbol{q}}_s$ become zero. Thus, $\tilde{\boldsymbol{q}}_s$ can be represented as $-\tilde{q}_s\hat{z}$, where $\tilde{q}_s=|{\tilde{\boldsymbol{q}}}_s|$. Consequently, the average direction of swimming of microorganisms can be determined by 
	\begin{equation*}
		<\boldsymbol{p_s}>=-M(G_s)\frac{\tilde{{\boldsymbol{q}}}_s}{\tilde{q}_s}=M_s\hat{z},
	\end{equation*}
	where $M_s=M(G_s)$.
	The concentration $n_s(z)$ satisfies
	\begin{equation}\label{32}
	\frac{dn_s}{dz}-V_cM_sn_s=0,
	\end{equation}
	which is accompanied by the cell conservation relation
	\begin{equation}\label{33}
	\int_0^1n_s(z)dz=1.
	\end{equation}
	
	A boundary value problem (BVP) is represented by equations (\ref{31}) through (\ref{33}), and is solved numerically by the shooting technique~\cite{15panda2016,41rajput2023}.
	
	\begin{figure*}
		\includegraphics[width=15cm]{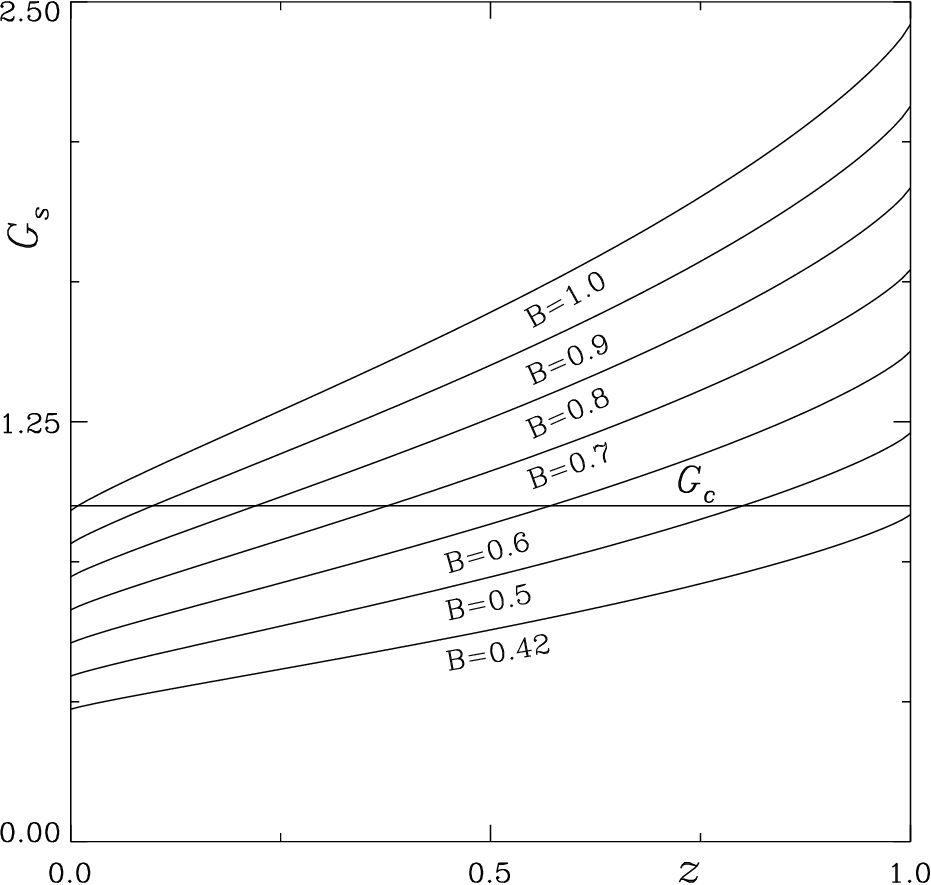}
		\caption{\label{fig2} Alteration in the total intensity for a homogeneous suspension at different values of $B$, while $k=0.5$, and $\omega=0.7$.}
	\end{figure*}	
	
	Assuming that the microorganisms under consideration are phototactic and resemble Chlamydomonas~\cite{34hill1989,15panda2016,7ghorai2010}. We adopt the same parameter values used in previous studies to enable straightforward comparison. These parameter values are listed in Table 1 in Panda $et$ $al$.~\cite{8panda2020}, and the approach for calculating the required radiation properties follows the methodology presented by Panda $et$ $al.$~\cite{15panda2016}. For a suspension with a depth of 0.5cm, the optical depth $\tau\in[0.25,1]$, and $\omega\in[0,1]$~\cite{33daniel1979}. Non-dimensional swimming speed is estimated to be $V_c=10$ and $V_c=20$ for suspensions with depths of 0.5cm and 1.0cm, respectively (refer to Table 1 in Panda $et$ $al$.~\cite{8panda2020}). The intensity of diffuse irradiation $B\in[0,1]$, depending on the cloud cover.	
	
	Consider the phototaxis function $M(G)$, whose mathematical form is given by
	
	\begin{equation}\label{34}
	M(G)=0.8\sin\left[\frac{3\pi}{2}\aleph(G)\right]-0.1\sin\left[\frac{\pi}{2}\aleph(G)\right],	
	\end{equation}
	
	where, $\aleph(G)=\frac{G}{2.5}\exp\left[0.32(2.5-G)\right]$ with $G_c=1$.
	
	\begin{figure*}[!htbp]
		\includegraphics[width=15cm]{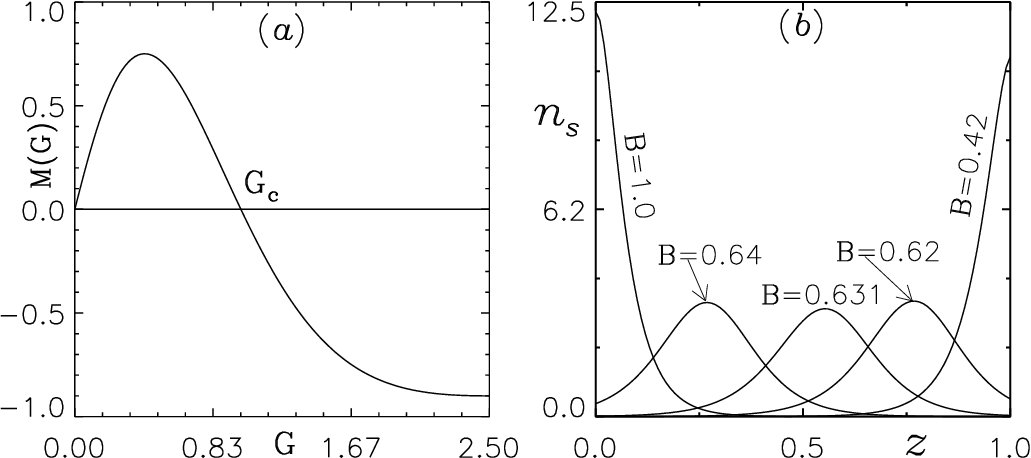}
		\caption{\label{fig3} (a) Photo-response curve for $G_c=1.0$, and (b) basic profiles of cell concentration for different sets of values of diffuse flux $B$, where $V_c=20$, $\kappa_H=0.5$, and $\omega=0.7$.}
	\end{figure*}
	
	Fig.~\ref{fig2} depicts the behaviour of total intensity $G_s$ throughout the suspension for different values of $B$ where other governing parameters such as $S_c,\kappa_H,\omega$ are kept fixed at 20,0.5 and 0.7, respectively. Here, $G_s$ monotonically decreases across the suspension for $0<B\leq 1$.
	Fig.~\ref{fig3} displays the photo-response curve for the critical intensity $G_c=1$ and basic profiles of cell concentration for $V_c=20$, $\omega=0.7$, and $\kappa_H=0.5$ at different values of $B$. For $0<B\leq 0.42$, the sublayer at steady state occurs at the upper boundary of the domain (i.e., at the top of the suspension). However, as $B$ increases beyond this range, the position of the sublayer relocates towards the lower boundary of the suspension (see Fig.~\ref{fig3}(b)).\par

	Now, we consider the case of purely scattering ($\omega=1$) in the suspension. Here, we assumed the taxis function 
	\begin{align}\label{35}
		M(G)=0.8\sin\left[\frac{3\pi}{2}\aleph(G)\right]-0.1\sin\left[\frac{\pi}{2}\aleph(G)\right], 	
	\end{align}
	
	where $\aleph(G)=\frac{G}{2.5}\exp[0.04(2.5-G)]$ with the critical intensity $G_c=1.55$.
	
	\begin{figure*}[!htbp]
		\includegraphics[width=15cm]{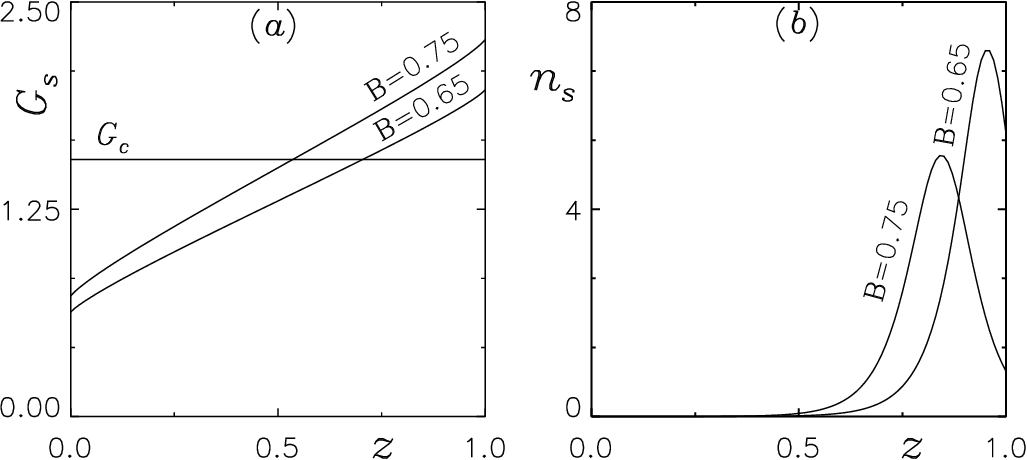}
		\caption{\label{fig4}(a) The variation in total intensity in a uniform suspension for two different values of $B$ and, (b) corresponding basic profiles of cell concentration in the suspension for the same values of $B$. The parameters $S_c$, $V_c$, $\kappa_H$, and $\omega$ are held constant at $20$, 20, $0.5$, and $1$, respectively.}
	\end{figure*}

	Fig.~\ref{fig4} shows the variation in $G_s$ in a uniform suspension $(n=1)$, where the suspension is highly scattering, where the other governing parameters $V_c=20$, $\kappa_H=1$ are held constant. The figure demonstrates that $G_s$ across the suspension decreases monotonically in highly scattering suspension, as in the previous case of weak scattering suspension, which is not seen in highly scattering algal suspension generally. Consequently, the critical intensity $G_c$ occurs at a single location in the suspension and an unimodal basic state is observed for highly scattering suspensions.
	
	
	\section{THE LINEAR STABILITY ANALYSIS }
	The linear perturbation theory is used for linear stability analysis (LSA), where a small disturbance (perturbation) of amplitude $\tilde{\epsilon}$ (where $0<\tilde{\epsilon}\ll1$) is made to the static state as follows

	\begin{align}\label{36}
		[\boldsymbol{U},n,L,<\boldsymbol{p}>]=[0,n_s,L_s,<\boldsymbol{p}_s>]+\tilde{\epsilon} [\boldsymbol{U}_1,n_1,L_1,<\boldsymbol{p}_1>]
		+O(\tilde{\epsilon}^2)
		.  
	\end{align}
	
	The linearization of Eqs. (\ref{15})-(\ref{17}) is conducted by using the perturbed quantities and gathering terms of order $O(\tilde{\epsilon})$. The resultant linearized equations can be represented as
	\begin{equation}\label{37}
	\boldsymbol{\nabla}\cdot \boldsymbol{U}_1=0,
	\end{equation}
	where  $\boldsymbol{U}_1=(U_1,V_1,W_1)$.
	\begin{equation}\label{38}
	S_{c}^{-1}\left[\frac{\partial \boldsymbol{U_1}}{\partial t}\right]=-\boldsymbol{\nabla} P-Rn_1\hat{\boldsymbol{z}}+\nabla^{2}\boldsymbol{U}_1,
	\end{equation}
	\begin{equation}\label{39}
	\frac{\partial{n_1}}{\partial{t}}+V_c\boldsymbol{\nabla}\cdot[<\boldsymbol{p_s}>n_1+<\boldsymbol{p_1}>n_s]+W_1\frac{dn_s}{dz}=\boldsymbol{\nabla}^2n_1.
	\end{equation}
	If $G=G_s+\tilde{\epsilon} G_1+O(\tilde{\epsilon}^2)=G_s^d+\tilde{\epsilon} G_1^d+O(\tilde{\epsilon}^2)$, then $G_1$ is given by
	\begin{equation}\label{40}
	G_1=G_1^d=\int_{4\pi}L_1^d(\boldsymbol{x},\boldsymbol{s})d\Omega.
	\end{equation}
	Similarly  
	\begin{equation}\label{41}
	\tilde{\boldsymbol{q}}_1=\tilde{\boldsymbol{q}}_1^d=\int_{4\pi}L_1^d(\boldsymbol{x},\boldsymbol{s})\boldsymbol{s}d\Omega.
	\end{equation}
	The perturbed
	\begin{equation}\label{42}
	<\boldsymbol{p_1}>=G_1\frac{dM_s}{dG}\hat{\boldsymbol{k}}-M_s\frac{\boldsymbol{\tilde{q}}_1^{H}}{\tilde{q}_s},
	\end{equation}
	where $\tilde{\boldsymbol{q}}_1^H$ is the horizontal component of  $\tilde{\boldsymbol{q}}_1$.
	Now, we substitute the value of $<p_1>$ in equation (\ref{39}) and simplify, obtain
	\begin{widetext}
		\begin{equation}\label{43}
		\frac{\partial{n_1}}{\partial{t}}+V_c\frac{\partial}{\partial z}\left[M_sn_1+n_sG_1\frac{dM_s}{dG}\right]-V_cn_s\frac{M_s}{{\tilde{q}_s}}\left[\frac{\partial \tilde{q}_1^x}{\partial x}+\frac{\partial \tilde{q}_1^y}{\partial y}\right]
		+W_1\frac{dn_s}{dz}-\nabla^2n_1=0.
		\end{equation}
	\end{widetext}
	To simplify the system of equations, the curl in Eqs.~(\ref{38}) is applied twice, and only the vertical component of the resulting equation is retained. As a result, $P$ and the horizontal component of $\boldsymbol{U}_1$ is eliminated, which reduces Eqs.~(\ref{37})-(\ref{39}) to two equations describing $W_1$ and $n_1$. These equations can be solved for normal modes
	
	\begin{equation}\label{44}
	[W_1,n_1]=[\tilde{W}(z),\tilde{\vartheta}(z)]\exp{[\sigma t+i(l_1x+l_2y)]}.  
	\end{equation}
	
	The diffuse perturbed intensity $L_1^d$ can be governed by the equation
	\begin{equation}\label{45}
	\alpha\frac{\partial L_1^d}{\partial x}+\beta\frac{\partial L_1^d}{\partial y}+\gamma\frac{\partial L_1^d}{\partial z}+\kappa_H( n_sL_1^d+n_1L_s^d)=\frac{\omega\kappa_H}{4\pi}(n_sG_1^d+G_s^dn_1),
	\end{equation}
	with the boundary conditions
	\begin{subequations}
		\begin{equation}\label{46a}
		L_1^d(x, y, z=1, \alpha,\beta, \gamma) =0,~~ where~~(\pi/2\leq\vartheta\leq\pi,~~0\leq\varphi\leq 2\pi), 
		\end{equation}
		\begin{equation}\label{46b}
		L_1^d(x, y, z=0, \alpha,\beta, \gamma) =0,~~where~~ (0\leq\vartheta\leq\pi/2,~~0\leq\varphi\leq 2\pi). 
		\end{equation}
	\end{subequations}
	
	The expression of Eq.~(\ref{45}) implies the following form for $L_1^d$	\begin{equation*}
		L_1^d=\Psi^d(z,\alpha,\beta,\gamma)\exp{(\sigma t+i(l_1x+l_2y))}. 
	\end{equation*}
	From Eq.~(\ref{40}), we get
	\begin{widetext}
		
		\begin{equation*}
			G_1^d=\mathcal{G}^d(z)\exp{(\sigma t+i(l_1x+l_2y))}=\left[\int_{4\pi}\Psi^d(z,\alpha,\beta, \gamma)d\Omega\right]\exp{(\sigma t+i(l_1x+l_2y))},
		\end{equation*}
	\end{widetext}
	where
	\begin{equation}\label{47}
	\mathcal{G}^d(z)=\int_{4\pi}\Psi^d(z,\alpha,\beta,\gamma)d\Omega.
	\end{equation}
	Now $\Psi^d$ satisfies
	\begin{widetext}
		\begin{equation}\label{48}
		\frac{d\Psi^d}{dz}+\frac{[i(l_1\alpha+l_2\beta)+\kappa_H n_s]}{\gamma}\Psi^d=\frac{\omega\kappa_H}{4\pi\gamma}(n_s\mathcal{G}^d+G_s^d\Theta)-\frac{\kappa_H}{\gamma}L_s^d\Theta,
		\end{equation}
	\end{widetext}
	with the boundary conditions
	\begin{subequations}
		\begin{equation}\label{49a}
		\varPsi^d( z=1,  \alpha,\beta, \gamma) =0,~~ (\pi/2\leq\vartheta\leq\pi,~~0\leq\varphi\leq 2\pi), 
		\end{equation}
		\begin{equation}\label{49b}
		\varPsi^d( z=0,\alpha,\beta, \gamma) =0,~~ (0\leq\vartheta\leq\pi/2,~~0\leq\varphi\leq 2\pi). 
		\end{equation}
	\end{subequations}
	From Eq.~(\ref{48}), we have
	\begin{equation*}
		\tilde{q}_1^H=[\tilde{q}_1^x,\tilde{q}_1^y]=[P(z),Q(z)]\exp{[\sigma t+i(l_1x+l_2y)]},
	\end{equation*}
	where
	\begin{align*}
		[P(z),Q(z)]=\int_0^{4\pi}(\alpha,\beta)\Psi^d(z,\alpha,\beta,\gamma) d\Omega.
	\end{align*}
	The governing equations become
	\begin{equation}\label{50}
	\left[\sigma S_c^{-1}+k^2-\frac{d^2}{dz^2}\right]\left[ \frac{d^2}{dz^2}-k^2\right]\tilde{W}-Rk^2\tilde{\Theta}=0,
	\end{equation}
	\begin{widetext}
		\begin{equation}\label{51}
		\left[\sigma+k^2-\frac{d^2}{dz^2}\right]\tilde{\Theta}+V_c\frac{d}{dz}\left[M_s\tilde{\Theta}+n_s\mathcal{G}^d\frac{dM_s}{dG}\right]\\
		-i\frac{V_cn_sM_s}{\bar{q}_s}(l_1P+l_2Q)+\frac{dn_s}{dz}\tilde{W}=0,
		\end{equation}
	\end{widetext}
	with the boundary conditions
	\begin{equation}\label{52}
	\frac{d\tilde{\Theta}}{dz}-V_cM_s\tilde{\Theta}-n_sV_c\mathcal{G}^d\frac{dM_s}{dG}=0, ~~~z=0,1.
	\end{equation}
	For rigid boundary
	\begin{equation}\label{53}
	\tilde{W}=\frac{d\tilde{W}}{dz}=0,~~~z=0
	\end{equation}
	and for the stress-free boundary,
	\begin{equation}\label{54}
	\tilde{W}=\frac{d^2\tilde{W}}{dz^2}=0,~~~z=1.
	\end{equation}
	
	After simplification, Eq.~(\ref{51}) becomes (using D=d/dz)
	\begin{equation}\label{55}
	\left[\sigma S_c^{-1}+k^2-D^2\right]\left[ D^2-k^2\right]\tilde{W}-Rk^2\tilde{\Theta}=0,
	\end{equation}
	\begin{equation}\label{56}
	\Gamma_0(z)+(\sigma+k^2+\Gamma_1(z))\tilde{\Theta}+V_cM_sD\tilde{\Theta}
	-D^2\tilde{\Theta}+Dn_s\tilde{W}=0, 
	\end{equation}
	where
	\begin{subequations}
		\begin{equation}\label{57a}
		\Gamma_0(z)=V_cD\left[n_s\mathcal{G}^d\frac{dM_s}{dG}\right]-i\frac{V_cn_sM_s}{\bar{q}_s}(l_1P+l_2Q),
		\end{equation}
		\begin{equation}\label{57b}
		\Gamma_1(z)=V_c\frac{dM_s}{dG}DG_s^d,
		\end{equation}
	\end{subequations}
	
	with the boundary conditions
	\begin{equation}\label{58}
	D\tilde{\Theta}-V_cM_s\tilde{\Theta}-n_sV_c\mathcal{G}^d\frac{dM_s}{dG}=0, ~~~z=0,1.
	\end{equation}
	For rigid boundary
	\begin{equation}\label{59}
	\tilde{W}=D\tilde{W}=0,~~~z=0
	\end{equation}
	and for the stress-free boundary,
	\begin{equation}\label{60}
	\tilde{W}=D^2\tilde{W}=0,~~~z=1.
	\end{equation}
	
	Here, terms $\Gamma_0$ and $\Gamma_1$ in Eq. (\ref{56}), are used to include the novel characteristics of the proposed model through the impacts of diffuse flux.
	
	\section{SOLUTION PROCEDURE}
	The set of Eqs.~(\ref{55}) and (\ref{56}) are solved by the Newton Raphson Kantorocich (NRK) method~\cite{19cash1980}. It is a numerical technique that involves iterative improvement of an initial guess until the desired accuracy is achieved. Using this method, one can determine the neutral (or marginal) stability curves. The neutral curve $R^{(n)}(k)$ has infinite branches and represents a unique solution for fixed parameter values. The most unstable (interesting) bioconvective solution corresponds to the branch where the Rayleigh number $R$ attains its minimum value, and the wavelength of instability is computed using $\lambda_c=2\pi/k_c$, where $k_c$ is the critical wave number. The bioconvective solutions are characterized by the presence of convection cells arranged vertically in the suspension, with the number of cells being indicated by the mode of the solution. 	
	\section{NUMERICAL RESULTS}
	To investigate how diffuse flux solely (without collimated flux) affects phototactic bioconvection, we conducted a parameter study by varying the value of $B$ from 0 to 1 while keeping other parameters $S_c$, $V_c$, $\kappa_H$, and $\omega$ are fixed constant. To manage a large number of parameters, we select a discrete set of constant values. We employ $S_c=20$ and varied ($V_c$) between 10 and 20, $\kappa_H$ from 0.5 to 1.0, $\omega$ from 0 to 1, and $B$ from 0 to 1, to cover a broad range of values. Since for all values of $\omega$, the total intensity shows the same behaviour, we select $\omega=0.7$. The value of parameter $B$ is selected to ensure that the location of maximum concentration (LMC) at the equilibrium state shifts from $z=1$ to $z=1/2$ height of the domain. The depth of the gravitationally unstable region below the LMC (DGUR) promoted convection, while the height of the gravitationally stable region (hereinafter referred to as HGSR) above the LMC inhibited convective fluid motion. The stability of the suspension is also affected by the steepness of the maximum concentration (SMC). As DGUR and SMC rise, the critical Rayleigh number ($R_c$) reduces, indicating that the suspension's stability decreases. Conversely, as HGSR increased, the $R_c$ increased, suggesting that the suspension's stability increased.
	\subsection{ $V_c=$ 20}		
	\begin{figure*}[!htbp]
		\includegraphics[width=16cm]{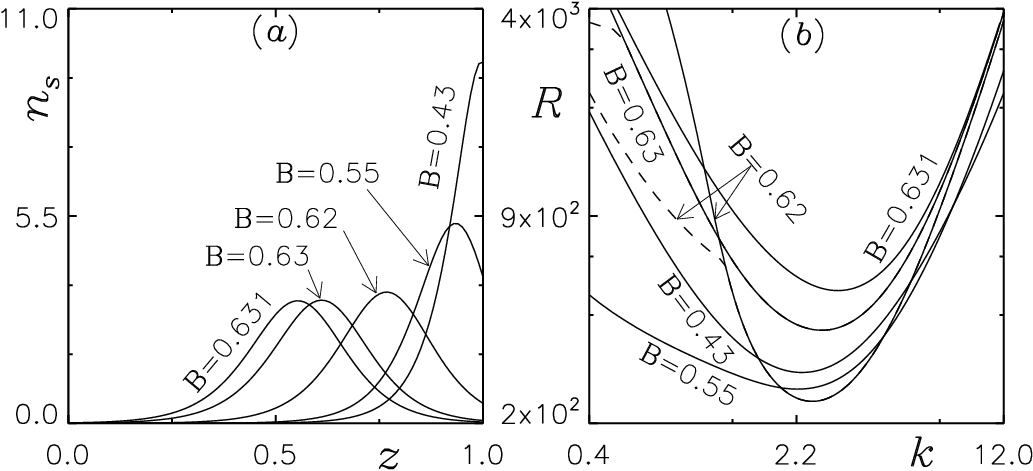}
		\caption{\label{fig5} (a) The basic profiles of cell concentration at steady state, (b) corresponding marginal stability curves. Here the upper surface is stress-free, and the parameter values $V_c=20$, $k_H=0.5$, and $\omega=0.7$.}
	\end{figure*}
\begin{figure*}[!htbp]
		\includegraphics[width=8.8cm]{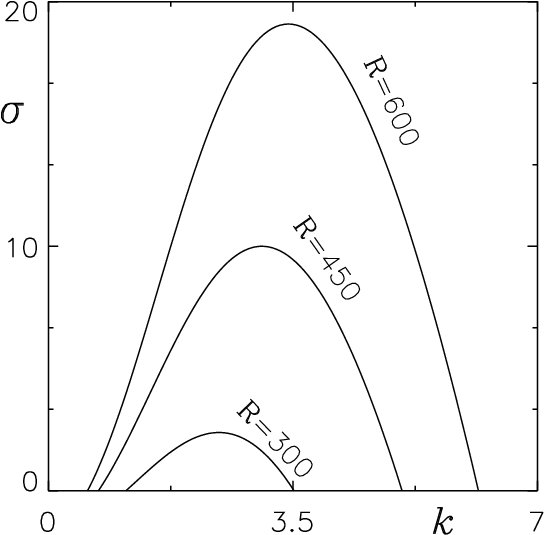}
		\caption{\label{fig6} Growth rate curves for $V_c=20$, $k_H=0.5$, $\omega=0.7$ and $B=0.43$ when $R>R_c$. Here the upper surface is stress-free.}
	\end{figure*}

		
	\textbf{(a)~~$\kappa_H= 0.5$}: Figure~\ref{fig5}(a) exhibits the cell concentration profiles at the basic state and their variations as $B$ increases while maintaining constant values for the governing parameters $V_c=20$, $\kappa_H=0.5$, and $\omega=0.7$. Figure~\ref{fig5}(b) showcases the marginal stability curves corresponding to the concentration profiles in Figure~\ref{fig5}(a). At $B=0.43$, the sublayer (maximum concentration) is observed at the top of the suspension, and the $R_c$ occurs at the stable branch of the corresponding neutral curve, which implies that the solution is stationary. Here, $R_c$ is
	approximately 264.94 and the critical wavenumber is 2.32. Thus, the
	critical wavelength is 2.71. Figure~\ref{fig6} displays the growth rate $\sigma$ for $R>R_c$ when $B=0.43$. The most unstable mode is the one with the highest growth rate. Thus, the pattern wavelength decreases with increasing $R$. For example, the wavenumbers of the most unstable mode are approximately 2.53 and 3.48 for $R=300$ and $R=600$ respectively [see Figure ~\ref{fig6}]. 		
	
	As $B$ rises, the sublayer shifts nearer to $z=1/2$ within the suspension, consequently increasing the HGSR. When $B=0.55$, the sublayer is situated around $z\approx 0.88$, further augmenting the HGSR, which impedes convection. Additionally, positive phototaxis in the lower region hampers bioconvection, while negative phototaxis in the upper region promotes it. The effect of negative phototaxis is more pronounced for $B=0.55$ compared to $B=0.43$, leading to a lower critical Rayleigh number. As $B$ reaches 0.62, the sublayer occurs at $z\approx 3/4$ within the suspension. In this scenario, a similar mechanism occurs, leading to a small $R_c$. Moreover, instability leads to the division of the stationary branch of the neutral curve into an oscillatory branch. Nevertheless, it's important to note that the most unstable mode is associated with the stationary branch, which implies that the solution is stationary.
	
	\begin{figure*}[!htbp]
		\includegraphics{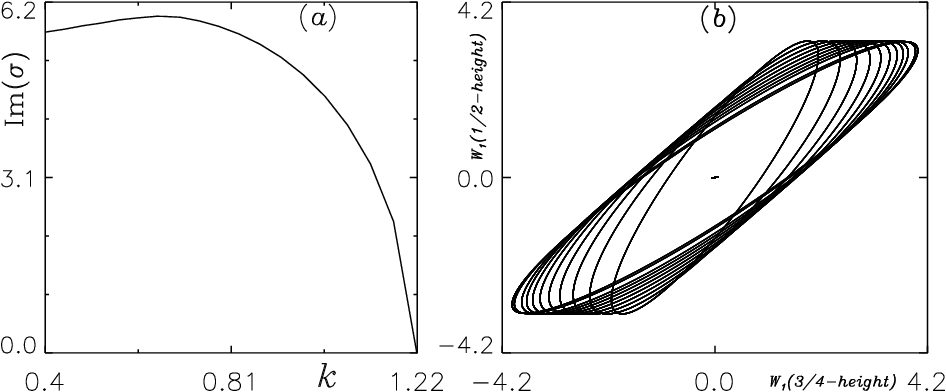}
		\caption{\label{fig7}(a) Variation of $Im(\sigma)$ vs wavenumber $k$ at the initiation of instability, and (b) limit cycle corresponding to the oscillatory branch of marginal stability curve for $B=0.62$. Here the upper surface is stress-free, and $V_c=20$, $k_H=0.5$, and $\omega=0.7$.}
	\end{figure*}
	
	In Figure \ref{fig7}(a), there is a connection between positive frequency and the wavenumber $k$ on the wavy branch, specifically when $B=0.62$. As the frequency gradually drops and reaches zero, the wavy motion changes into a still state at the start of the bioconvection process. This change is visible in Figure \ref{fig7}(a). Figure \ref{fig7}(b) shows a cool pattern: the motion gets smaller with each round, like a shrinking spiral. This means that the wavy flow becomes less and less bouncy over time, leading to a steady flow as the bouncing frequency gets closer to zero.

	As $B$ increases to 0.63, the sublayer forms around $z=0.6$ within the domain. Here, the HGSR increases and a higher $R_c$ occurs. The stable branch of the marginal stability curve divides into an oscillatory branch, and it's important to note that the most unstable solution corresponds with the stationary branch, which implies that the solution is stationary. The positive frequency shows the behaviour same as the case of $B=0.62$. The oscillatory branch transitions into a stationary branch as the positive frequency reaches zero. For $B=0.631$, the sublayer is situated at $z\approx1/2$ within the domain. In this scenario, the HGSR is high compared to the previous cases, resulting in the largest critical Rayleigh number. Furthermore, the oscillatory branch disappears and the solutions remain stationary.	
	\begin{figure*}[!htbp]
		\includegraphics[width=16cm]{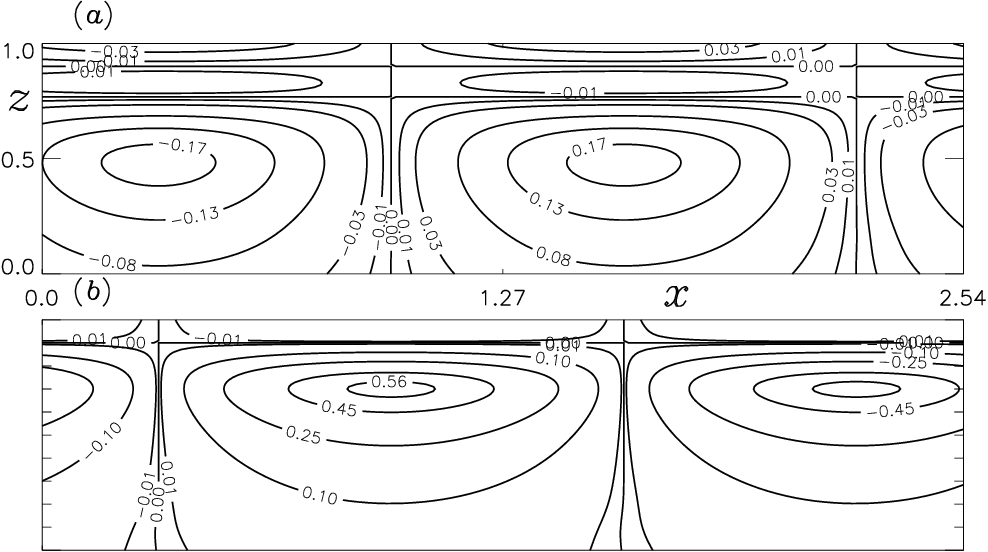}
		\caption{\label{fig8} The components of swimming orientation (a) horizontal component $<p_x>$, and (b) vertical component $<p_z>$. The parameter values $V_c=20$, $k_H=0.5$, $\omega=0.7$, and $B=0.62$.}
	\end{figure*}	
	\begin{figure*}[!htbp]
		\includegraphics[width=16cm]{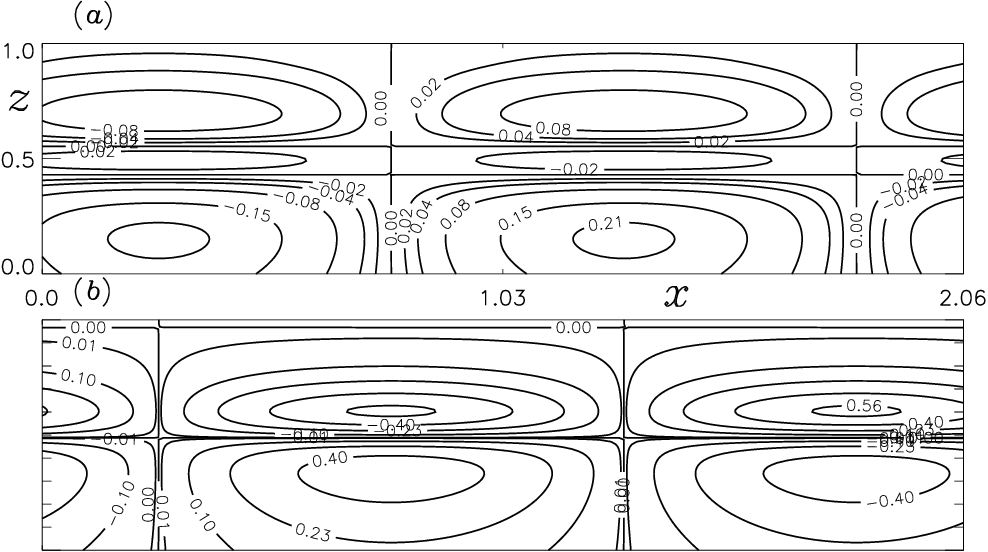}
		\caption{\label{fig9} The components of swimming orientation (a) horizontal component $<p_x>$, and (b) vertical component $<p_z>$. The parameter values $V_c=20$, $k_H=0.5$, $\omega=0.7$, and $B=0.631$ are kept constant.}
	\end{figure*}		
	
	Figure~\ref{fig8} illustrates a contour plot representing the components of the average swimming orientation $<\boldsymbol{p}>$ for a specific value of $B=0.62$. In this particular case, the contour plot provides a visual depiction of the typical swimming orientation components. The highest magnitudes observed for the horizontal component $<\boldsymbol{p}_x>$ and the vertical component $<\boldsymbol{p}_z>$ of the mean swimming orientation across the entire domain are approximately $1.7\times 10^{-1}$ and $5.6\times 10^{-1}$, respectively. It is noteworthy that $<\boldsymbol{p}_x>$ appears to be stronger in this case due to the uniformity of diffuse flux. The main reason behind this phenomenon is that diffuse flux contributes to both the horizontal and vertical components of swimming orientation, while collimated flux only affects the vertical component.		
	\begin{figure*}[!htbp]
		\includegraphics[width=16cm]{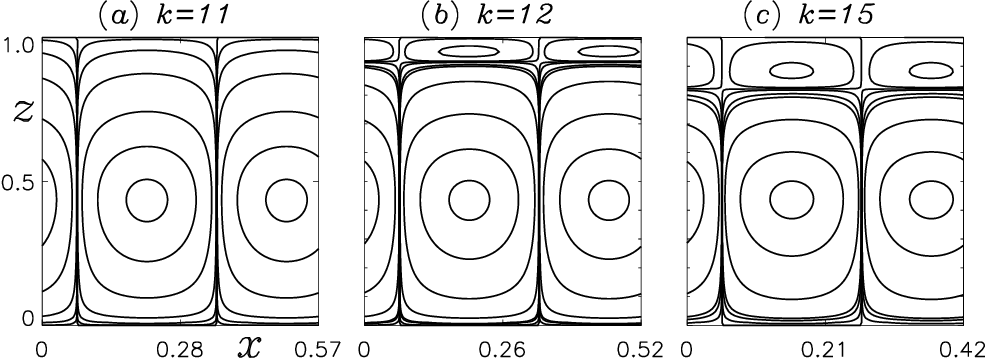}
		\caption{\label{fig10} The presence of mode 2 solutions during the bioconvective instability and the pattern of perturbed swimming velocity component $w_1$ are observed under constant parameter values of $V_c=20$, $\kappa_H=0.5$, $\omega=0.7$, and $B=0.631$. At the upper region of the suspension, a small convection cell becomes visible, and it enlarges as the wavenumber $k$ surpasses the critical value $k_c$. Here the upper surface is stress-free.}
	\end{figure*}	  		
	When $B=0.631$, the contour plots of $<\boldsymbol{p}_x>$ and $<\boldsymbol{p}_z>$ display in Figure~\ref{fig9}. It's worth noting that in this particular case, the critical wavelength is relatively small, as indicated in Table~\ref{tab2}. In this scenario, the maximum magnitudes of these components are approximately $2.1\times 10^{-1}$ and $5.6\times 10^{-1}$, respectively. Therefore, at the small wavelength, the diffuse flux contributes more to the mean horizontal orientation.

	The $R^{(1)}(k)$ branch displays a mode 2 bioconvective solution instead of mode 1 as diffuse flux rises. Specifically, as the wavenumber $k\geq k_c$, a smaller convection cell is observed to form at the top of the suspension and subsequently grows in height as $k$ increases more. This behaviour is consistently observed for $B=0.63$ and $B=0.631$. However, when $B=0.63$ and $B=0.631$, the entire branch is characterized by mode 2 behaviour, except for a specific range where the solution exhibits mode 1 behaviour. As a result, the bioconvective solution is denoted as a mode 2 solution (see Fig.~\ref{fig10}).
	
	
	\begin{figure*}[!htbp]
		\includegraphics[width=16cm]{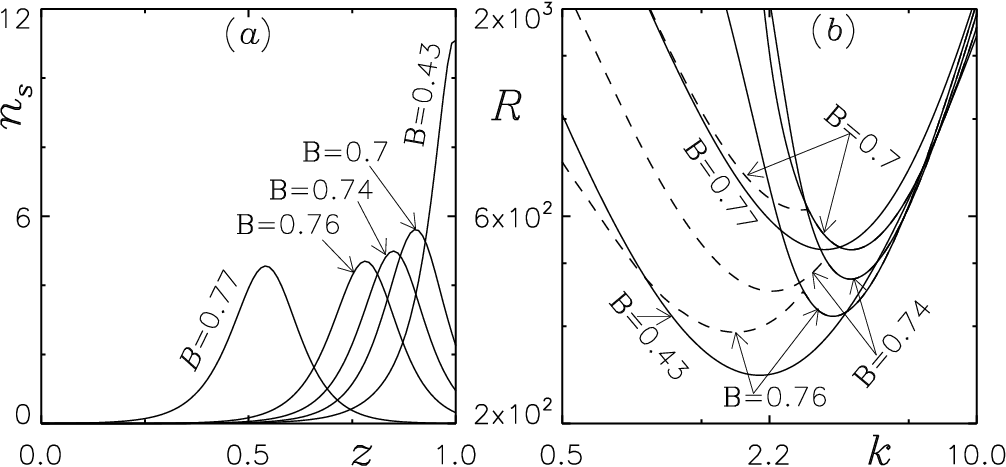}
		\caption{\label{fig11}(a) The basic profiles of cell concentration of concentration at steady state, (b) corresponding marginal stability curves. Here the upper surface is stress-free, and the parameter values $V_c=20$, $k_H=1$, and $\omega=0.7$ are kept constant.}
	\end{figure*}
	
	\textbf{(b)~~$\kappa_H=1$}: Figure~\ref{fig11} visually depicts the variations in the basic profiles of cell concentration alongside the neutral curves for different values of $B$, while maintaining the other parameters fixed at $V_c=20$, $\kappa_H=1$, and $\omega=0.7$. When $B=0.43$, the sublayer is situated at the top boundary of the suspension. The LSA predicts that the bioconvective solution remains stationary. As $B$ increases, the maximum cell concentration gradually shifts towards $z=1/2$ within the domain. At $B=0.7$, the sublayer moves to $z=0.9$, leading to an increase in the HGSR. This rise in HGSR hampers convective fluid motion, resulting in an increased critical Rayleigh number. Furthermore, an oscillatory branch emerges from the stationary branch of the neutral curve. However, it's important to note that the $R_c$ still resides on the stable branch, leading to a stable solution. When $B=0.74$, the sublayer observes at $z\approx 0.85$, and the HGSR increases. Furthermore, positive phototaxis in the lower region hinders the occurrence of convection, on the other hand, negative phototaxis in the upper region promotes it. The interplay between these factors results in the bioconvective solution becoming overstable. The bioconvective instability occurs at $k_c\approx2.28$ and $R_c\approx411.67$ with a non-zero $\sigma=\pm 17.20$. This type of bifurcation is referred to as periodic or Hopf bifurcation. Furthermore, the impact of negative phototaxis in the UGSR is stronger for $B=0.74$ compared to $B=0.7$, leading to a lower $R_c$.   	
	\begin{figure*}[!htbp]
		\includegraphics[width=14.2cm]{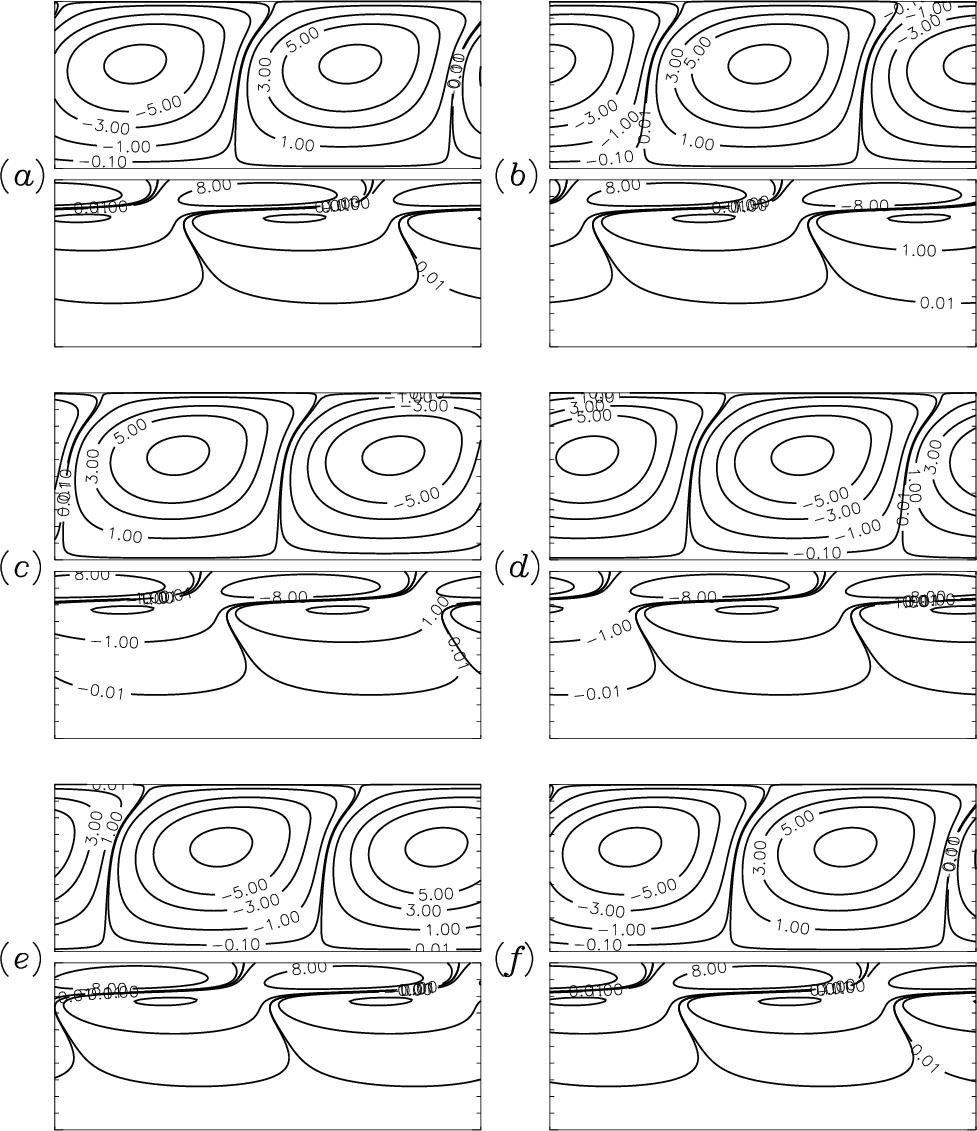}
		\caption{\label{fig12} The flow pattern of perturbed velocity component $w_1$ (top) and perturbed concentration component $n_1$ (bottom) for $B=0.74$ during a cycle of oscillation at initiation of instability. Here the upper surface is stress-free, and the other parameter values $V_c=20$, $k_H=1$, $\omega=0.7$ are kept constant.}
	\end{figure*}	
	The flow patterns of $w_1$ and $n_1$ throughout a single cycle of oscillation at the commencement of instability are represented in Figure~\ref{fig12}. The time intervals between figures (a) to (f) are consistent, and the period is measured at 0.365 units. Figure \ref{fig13}(a) demonstrates the fluctuation of the perturbed fluid velocity component $w_1$ over time. On the other hand, Figure \ref{fig13}(b) shows the corresponding limit cycle at approximately $k_c\approx 2.28$.		
	\begin{figure*}[!htbp]
		\includegraphics[width=16cm]{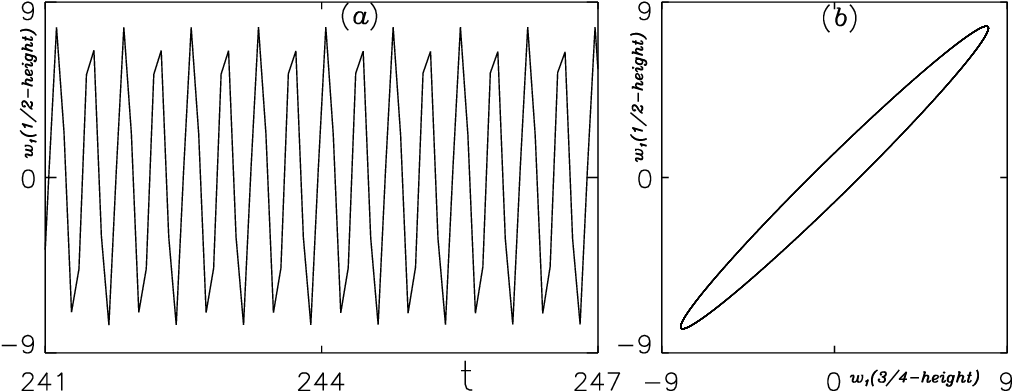}
		\caption{\label{fig13} (a) The perturbed ﬂuid velocity
			$w_1$ with respect to time and (b) limit cycle corresponding to the oscillatory branch of marginal stability curve for $B=0.74$, where $V_c=20$, $\kappa=1$, and $\omega=0.7$. Here the upper surface is stress-free. }
	\end{figure*}	
	When the magnitude of diffuse flux, $B=0.76$, the sublayer at the base state is observed around $z=3/4$ in the suspension. In this scenario, similar mechanisms are in play, resulting in a lower $R_c$ compared to the case of $B=0.74$. The oscillatory branch of the marginal stability curves indicates the presence of an overstable solution. Further increasing the magnitude of diffuse flux to 0.77 leads to the accumulation of cells at the steady state, occurring roughly at the mid-height of the suspension. It's worth noting that the HGSR is observed high in this case, resulting in the highest $R_c$. However, the nature of the bioconvective solution shifts from overstable to stable at $B=0.77$. In this context, the mode 2 instability exclusively occurs when $B=0.77$. A summary of the numerical findings is presented in Table~\ref{tab2}.
	
	\begin{table}[!bt]
		\caption{\label{tab2}The numerical results for $V_c=20$, with an increase in $B$, while keeping other parameters constant.}
		\begin{ruledtabular}
			\begin{tabular}{cccccccc}
				$V_c$ & $\kappa_H$ & $\omega$ & $B$ & $\lambda_c$ & $R_c$ & $Im(\sigma)$ & mode\\
				\hline
				\vspace{-0.2cm}\\			
				20 & 0.5 & 0.7 & 0.43 & 2.71 & 264.94 & 0 &  1\\
				20 & 0.5 & 0.7 & 0.55 & 2.85 & 233.50 & 0 &  1\\
				20 & 0.5 & 0.7 & 0.62\footnotemark[1] & 2.55 & 212.78 & 0 & 1\\
				20 & 0.5 & 0.7 & 0.63\footnotemark[1] & 2.35 & 365.23 & 0 & 2\\
				20 & 0.5 & 0.7 & 0.631 & 2.05 & 493.65 & 0 & 2\\
				
				20 & 1 & 0.7 & 0.43 & 3.01 & 260.20 & 0 &  1\\
				20 & 1 & 0.7 & 0.7\footnotemark[1] & 4.06 & 516.15 & 0 &  1\\
				20 & 1 & 0.7 & 0.74 & 2.76\footnotemark[2] & 411.67\footnotemark[2] & 11.20 &  1\\
				20 & 1 & 0.7 & 0.76 & 3.60\footnotemark[2] & 329.16\footnotemark[2] & 12.67 & 1\\
				20 & 1 & 0.7 & 0.77  & 1.87 & 516.51 & 0 &  2\\
				
			\end{tabular}
		\end{ruledtabular}
		\footnotetext[1]{The result indicates the existence of the oscillatory branch of the marginal stability curve.}
		\footnotetext[2]{The result indicates that the most unstable solution occurs on the oscillatory branch.}
	\end{table}
	
	
	\subsubsection{$V_c=10$ $and$ $15$}
	The impact of diffuse flux at the instability demonstrates similar qualitative behaviour for $V_c = 10$ and 15, as observed for $V_c = 20$. Notably, in this section, there are no instances of overstable solutions for any of the governing parameters. The numerical results can be found in Table~\ref{tab3}.
	
	\begin{table}[h]
		\caption{\label{tab3}The numerical results for $V_c=10$ and $V_c=15$, with an increase in $B$, while keeping other parameters constant.}
		\begin{ruledtabular}
			\begin{tabular}{cccccccc}
				$V_c$ & $\kappa_H$ & $\omega$ & $B$ & $\lambda_c$ & $R_c$ & $Im(\sigma)$ & mode\\
				\hline
				\vspace{-0.2cm}\\
				10 & 0.5 & 0.7 & 0.42 & 4.44 & 135.52 & 0 & 1 \\
				10 & 0.5 & 0.7 & 0.54 & 4.05 & 152.36 & 0 & 1\\
				10 & 0.5 & 0.7 & 0.59 & 3.77 & 201.41 & 0 & 1\\
				10 & 0.5 & 0.7 & 0.62 & 2.69 & 557.75 & 0 & 2\\
				10 & 0.5 & 0.7 & 0.63 & 2.06 & 1038.40 & 0 & 2\\
				10 & 1 & 0.7 & 0.42 & 3.57 & 168.72 & 0 &  1\\
				10 & 1 & 0.7 & 0.67 & 2.63 & 231.21 & 0 &  1\\
				10 & 1 & 0.7 & 0.72\footnotemark[1] & 2.51 & 255.32 & 0 &  1\\
				10 & 1 & 0.7 & 0.76 & 2.35 & 470.04 & 0 &  2\\
				10 & 1 & 0.7 & 0.77 & 2.01 & 761.40 & 0 &  2\\
				
				15 & 0.5 & 0.7 & 0.43 & 3.37 & 183.14 & 0 & 1 \\
				15 & 0.5 & 0.7 & 0.57 & 3.06 & 179.62 & 0 & 1\\
				15 & 0.5 & 0.7 & 0.61\footnotemark[1] & 2.90 & 205.98 & 0 & 1\\
				15 & 0.5 & 0.7 & 0.63\footnotemark[1] & 2.57 & 365.70 & 0 & 2\\
				15 & 0.5 & 0.7 & 0.632 & 1.85 & 803.13 & 0 & 2\\
				15 & 1 & 0.7 & 0.43 & 2.93  & 232.67  & 0 & 1\\
				15 & 1 & 0.7 & 0.7\footnotemark[1] & 1.94  & 331.85 & 0 & 1\\
				15 & 1 & 0.7 & 0.75\footnotemark[1] & 2.06 & 295.08 & 0 & 1\\
				15 & 1 & 0.7 & 0.76\footnotemark[1] & 2.15 & 301.97 & 0 & 1\\
				15 & 1 & 0.7 & 0.77 & 1.88 & 570.49 & 0 & 2\\
				
			\end{tabular}
		\end{ruledtabular}
		\footnotetext[1]{The result indicates the existence of the oscillatory branch of the marginal stability curve.}
	\end{table}
	
	
	\subsubsection{Effect of cell swimming speed ($V_c$)}
	
\begin{figure*}[!htbp]
		\includegraphics[width=16cm]{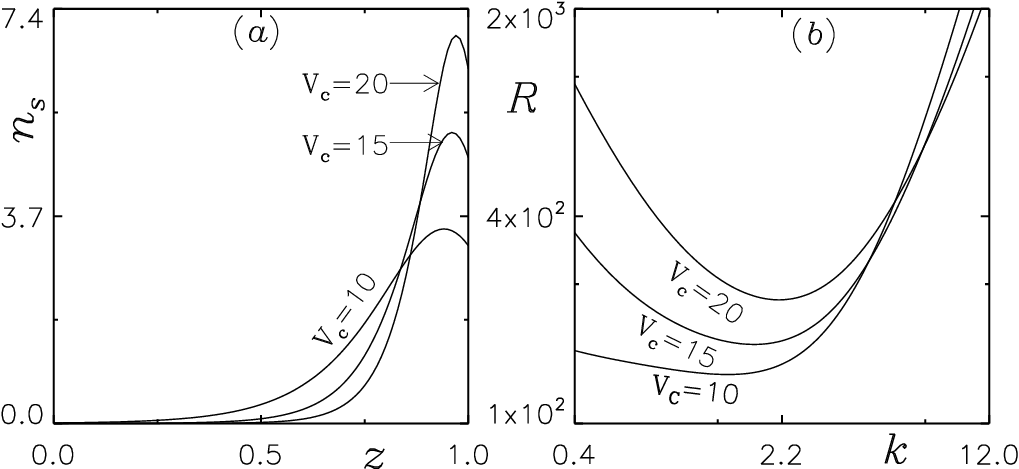}
		\caption{\label{fig14}(a) The basic profiles of cell concentration at steady state, (b) corresponding marginal stability curves. Here the upper surface is stress-free and the parameter values $V_c=20$, $k_H=0.5$, $\omega=0.7$, and $B=0.5$.}
	\end{figure*}

	In this section, we investigate how the swimming speed of cells $V_c$, affects the onset of bio-convection by varying the values of $B$ (0.5, 0.6, and 0.625) while keeping the governing parameters ($S_c=20$, $\kappa_H=0.5$, and $\omega=0.7$) constant. In the first case ($B=0.5$), when $V_c=10$, the sublayer occurs at approximately $z=0.94$. As $V_c$ increases, it shifts towards the top of the suspension, and SMC increases. However, at higher $V_c$, the cells encounter greater resistance due to positive phototaxis. Here, the second effect due to positive phototaxis dominates the former effect due to SMC. As a result, higher $R_c$ is observed for higher $V_c$ (see Fig.~\ref{fig14}).
	
\begin{figure*}[!htbp]
		\includegraphics[width=16cm]{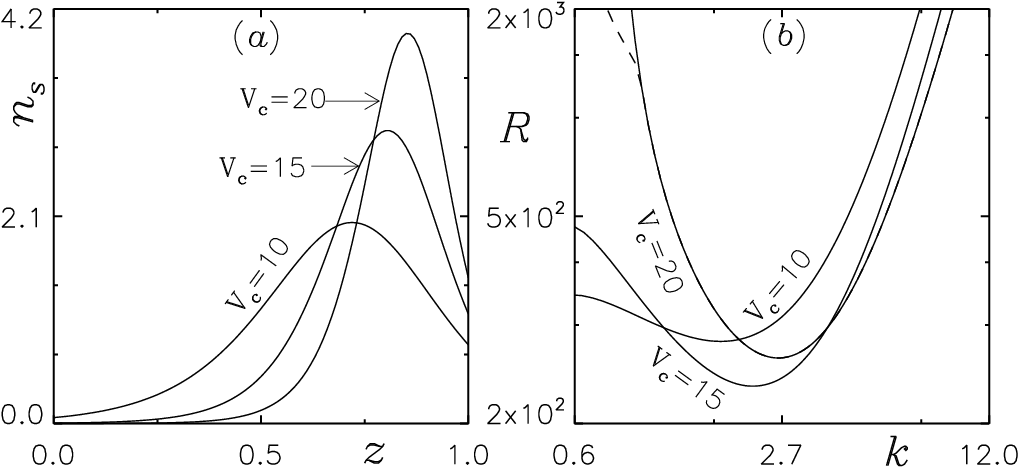}
		\caption{\label{fig15}(a) The basic profiles of cell concentration at steady state, (b) corresponding marginal stability curves. Here the upper surface is stress-free, and the parameter values $V_c=20$, $k_H=0.5$, $\omega=0.7$, and $B=0.6$.}
	\end{figure*}

	For the second case ($B=0.6$), when $V_c=10$, the sublayer arises around three-quarters of the domain. As $V_c$ rises, the sublayer relocates towards the top of the suspension, which leads to higher SMC and higher DGUR, both of which support convection. However, positive phototaxis becomes more effective at faster $V_c$, inhibiting convection. In this case, the latter effect dominates, resulting in a lower $R_c$ at $V_c=20$ compared to $V_c=15$. The $R_c$ is smaller at $V_c=15$ than at $V_c=10$ because the former effect dominates over the latter (see Fig.~\ref{fig15}). When $V_c=20$, it is also evident that a non-stationary branch emerges from the non-oscillatory segment of the neutral curve. However, the primary non-oscillatory solution persists as the most unstable bioconvective solution on this branch.
	
	\begin{figure*}[!htbp]
		\includegraphics[width=16cm]{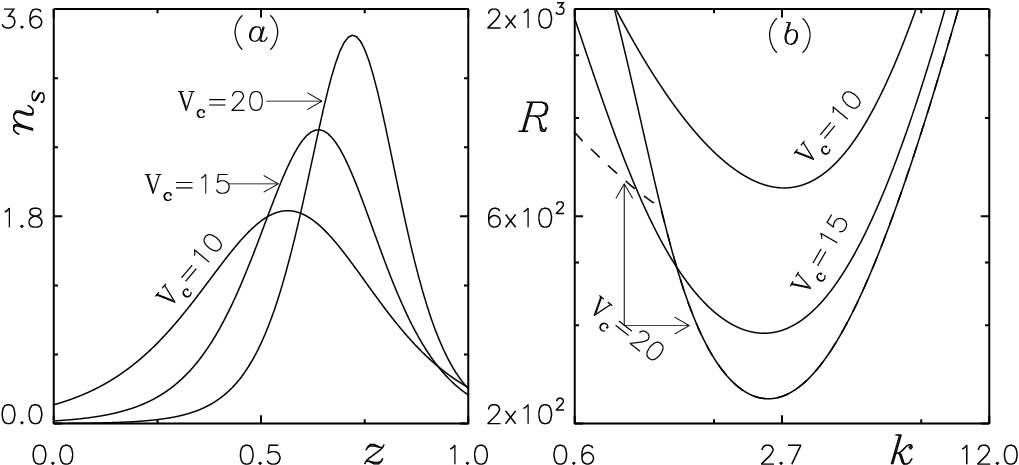}	
		\caption{\label{fig16}(a) The basic profiles of cell concentration at steady state, (b) corresponding marginal stability curves. Here the upper surface is stress-free, and the parameter values $V_c=20$, $k_H=0.5$, $\omega=0.7$, and $B=0.625$.}
	\end{figure*}	
	
	In the third case ($B=0.625$), the sublayer appears near $z=1/2$ height in the suspension, when $V_c=10$. As $V_c$ increases, the sublayer towards $z=3/4$ height of the suspension. Therefore, the SMC and HGSR both increase, supporting convection. However, positive phototaxis becomes stronger as $V_c$ increases from 10 to 20, inhibiting convection. In this case, the former effect prevails, resulting in a lower $R_c$ at $V_c=20$ compared to $V_c=10$ (see Fig.~\ref{fig16}). For $V_c=20$, a new oscillatory branch is split from the non-oscillatory branch of the neutral curve. However, the most unstable bioconvective solution remains on the non-oscillatory branch, leading to a stationary solution.
	The investigation of bioconvection instability with $\kappa_H=1$ reveals that the impact of cell swimming speed is comparable to that observed with $\kappa_H=0.5$. The bioconvective solutions obtained for $\kappa_H=1$ exhibit quantitative similarities to the solutions obtained for $\kappa_H=0.5$. 
	
	Detailed results and findings regarding bioconvection instability in this scenario can be found in Table~\ref{tab4}.
	
	\begin{table}[!htbp]
		\caption{\label{tab4}The numerical values of the bioconvective solutions for showing the effect of variation in cell swimming speed are shown in the table while keeping other parameters constant.}
		\begin{ruledtabular}
			\begin{tabular}{ c c c c c c }
				$B$  & $\kappa_H$ & $\omega$ & $V_c$ & $\lambda_c$ & $R_c$ \\
				\hline
				\vspace{-0.2cm}\\
				
				0.5 & 0.5 & 0.7 & 10 & 4.53 & 142.38  \\
				
				0.5 & 0.5 & 0.7 & 15 & 3.65 & 176.99 \\
				
				0.5 & 0.5 & 0.7 & 20 & 2.98 & 244.3 \\
				
				0.6 & 0.5 & 0.7 & 10 & 3.65 & 239.72 \\
				
				0.6 & 0.5 & 0.7 & 15 & 2.9 & 185.66 \\
				
				0.6 & 0.5 & 0.7 & 20\footnotemark[1] & 2.41 & 218.09  \\
				
				0.625 & 0.5 & 0.7 & 10 & 2.31 & 761.29 \\
				
				0.625 & 0.5 & 0.7 & 15 & 2.67 & 333.89 \\
				
				0.625 & 0.5 & 0.7 & 20\footnotemark[1] & 2.57 & 230.25 \\
				
				0.6 & 1 & 0.7 & 10 & 3.33 & 205.24  \\
				
				0.6 & 1 & 0.7 & 15 & 2.63 & 294.26 \\
				
				0.6 & 1 & 0.7 & 20 & 2.22 & 448.35  \\
				
				0.7 & 1 & 0.7 & 10 & 2.51 & 242.06\\
				
				0.7 & 1 & 0.7 & 15\footnotemark[1] & 1.94 & 331.85\\
				
				0.7 & 1 & 0.7 & 20\footnotemark[1] & 1.55 & 516.151 \\
				
				0.765 & 1 & 0.7 & 10 & 2.22 & 586.55 \\
				
				0.765 & 1 & 0.7 & 15 & 2.21 & 337 \\
				
				0.765 & 1 & 0.7 & 20\footnotemark[1] & 1.92 & 338.87 \\
				
			\end{tabular}
		\end{ruledtabular}
		\footnotetext[1]{The result indicates the existence of the oscillatory branch of the marginal stability curve.}
		\footnotetext[2]{The result indicates that the most unstable solution occurs on the oscillatory branch.}
	\end{table}
	
	\section{Comparison with similar bioconvection model}
	
	\begin{figure*}[!htbp]
		\includegraphics[width=16cm]{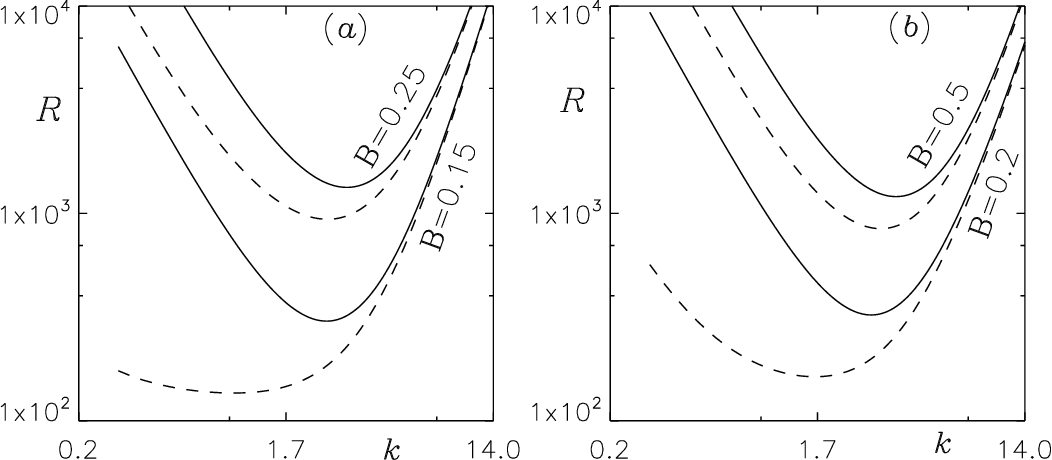}	
		\caption{\label{fig17}Comparison of marginal stability curves between proposed model (solid lines) and the model of both diffuse and collimated flux (dashed lines), (a) and (b) correspond to $B=0.27$ and $B=0.5$, respectively. Fixed parameter values are $V_c=10$, $\kappa_H=1$ and $\omega=0.4$.}
	\end{figure*}
	
	Here, we evaluate our model against a similar model on phototactic bioconvection in an algal suspension that is exposed to both diffuse and collimated flux~$\cite{15panda2016}$. To align the basic states of both models at the same location within the suspension, we transfer the collimated part of the basic total intensity from the model~$\cite{15panda2016}$ to our model and adjust it accordingly.
	
	Figure~$\ref{fig17}(a)$ presents the comparison of marginal stability curves between the two models for $B=0.15$ and $B=0.25$, with other governing parameters $V_c=10$, $\kappa_H=0.5$, and $\omega=0.4$ fixed. The marginal stability curves indicate agreement at small wavelengths but disagreement at small wavenumbers.
	
	For another set of parameters, Figure~$\ref{fig17}(b)$ illustrates the comparison of marginal stability curves for $B=0.2$ and $B=0.5$, with fixed parameters $V_c=10$, $\kappa_H=1$, and $\omega=0.4$. It is noteworthy that in the proposed model, light scattering is more effective due to the presence of only diffuse flux. This leads to a discrepancy between the two models at large wavelengths.

	
	\section{CONCLUSIONS}
	
	In this study, we introduced a new model for phototaxis that considers the impact of scattered or diffuse light flux in the absence of direct light on a suspension of algae with isotropic scattering. We used linear stability analysis to explore how this type of suspension behaves.
	
	Without direct collimated flux, the diffuse flux spreads uniformly, creating consistent light levels throughout the algae suspension. This means there is less variation in light across the suspension, and algae do not need to actively move to find the best light for photosynthesis. Unlike situations with focused light, there is no self-shading where algae block light from each other. The light intensity decreases uniformly across the depth of the suspension due to some light being absorbed by algae and some being scattered. This decrease affects the formation of patterns in the suspension, making them develop more slowly, if at all.
	
	The analysis showed that the suspension can have both stationary (non-oscillatory) and oscillatory (non-stationary) solutions when bioconvection starts. However, oscillatory solutions only happen when the swimming speed of cells is higher and there is more light absorption. A higher diffuse flux level makes the suspension more stable. In simpler terms, more scattered radiation creates a more stable environment for these microorganisms. The speed at which cells swim also affects stability, with lower speeds improving stability in low-light conditions but decreasing stability when there is plenty of light. The wavelength of the initial disturbance decreases with more diffuse light, supporting experimental results by Williams and Bees~\cite{4williams2011}.
	
	It is also observed that diffuse light also affects cell swimming orientation, especially at smaller pattern wavelengths, making their horizontal movement more significant, similar to gyrotaxis in certain cases.
	
	We compared our model with one that considers both diffuse and focused light. They agree well at short wavelengths but differ at longer wavelengths due to diffuse light.
	
	To validate our model, comparing it with experimental data from purely phototactic algal suspensions is crucial. However, finding such data is very challenging because most algae exhibit other behaviours alongside phototaxis. Therefore, identifying microorganisms primarily showing phototaxis is necessary. Despite this limitation, our model can simulate phototactic bioconvection in anisotropic algal suspensions using the appropriate phase function.

	\section*{Conflict of interest}
	The author affirms that there are no conflicts of interest.
	\section*{Data Availability}
	The article provides all the necessary supporting data for the plots and other findings within its contents.
	\nocite{*}
	\bibliography{diffuse_free}
	
\end{document}